\documentclass[11pt]{article}
\usepackage{mathtools}
\usepackage{relsize}
\usepackage{mathrsfs}
\usepackage{amsfonts,amsmath, amssymb,amsthm,amscd}
\usepackage[height=9.6in,width=5.95in]{geometry}
\usepackage{verbatim}
\usepackage{hyperref}
\usepackage{latexsym}
\usepackage{indentfirst}

\include{bibtex}

\newtheorem{theorem}{Theorem}[section]

\newtheorem{lemma}[theorem]{Lemma}
\newtheorem{proposition}[theorem]{Proposition}
\newtheorem{conjecture}[theorem]{Conjecture}
\newtheorem*{theorem*}{Theorem}

\theoremstyle{remark}
\newtheorem{definition}[theorem]{Definition}
\newcommand{\KAirytwo}{K_{2}}
\DeclareMathOperator{\Ai}{Ai}

\newcommand{\C}{\mathbb{C}}
\newcommand{\A}{\mathcal{A}}
\renewcommand{\L}{\mathcal{L}}
\newcommand{\1}{\mathbf{1}}

\newcommand{\E}{\mathbb{E}}
\newcommand{\N}{\mathbb{N}}
\newcommand{\Q}{\mathbb{Q}}
\newcommand{\Z}{\mathbb{Z}}
\newcommand{\R}{\mathbb{R}}
\renewcommand{\P}{\mathbb{P}}

\newcommand{\eps}{\varepsilon}
\def\P{\mathbb{P}}
\def\E{\mathbb{E}}

\newcommand{\KAirytwoext}{K^{{\rm ext}}_{2}}

\begin {document}
\title{Ergodicity of the Airy line ensemble}
\author{Ivan Corwin and Xin Sun}
\date{}
\maketitle

\begin{abstract}
In this paper, we establish the ergodicity of the Airy line ensemble with respect to
        horizontal shifts. This shows that it is the only candidate for Conjecture 3.2 in \cite{corwin2014brownian}, regarding the classification of ergodic line ensembles satisfying a certain Brownian Gibbs property after a parabolic shift.
\begin{flushleft}
\textbf{Keywords:}  the Airy line ensemble, ergodicity, Gibbs measure, extremal
\end{flushleft}
\end{abstract}


\section{Introduction} \label{sec:Introduction}

The \emph{Airy$_2  $ process} was introduced in  \cite{prahofer2002scale} and describes fluctuations in  random matrices, random surfaces and KPZ growth models. For instance, the Airy$_2 $ process describes the scaling limit of the largest eigenvalue in GUE Dyson's Brownian motion as the number of eigenvalues goes to infinity. For more information and background, see \cite{prahofer2002scale,johansson2003discrete,corwin2014brownian,corwin2013continuum,quastel2013airy} and reference therein.

One illuminating way to view the Airy$_2  $ process is as the top line of \emph{the Airy line ensemble} whose finite dimensional distributions are described by a determinantal process with the \emph{extended Airy$_2  $  kernel} as its correlation kernel. This determinantal process was introduced in \cite{prahofer2002scale} and the  existence of  a continuous version was established in \cite{corwin2014brownian}. In the Dyson's Brownian motion context, the  Airy line ensemble describes the limiting measure focusing on the top collection of evolving eigenvalues.

 The Airy$_2$ process is stationary with respect to horizontal shifts and Equation (5.15) in \cite{prahofer2002scale} showed that it satisfies the strong mixing condition and hence is ergodic. In this paper, we extend this result to the Airy line ensemble. This shows that it is the only candidate for Conjecture 3.2 in \cite{corwin2014brownian}, regarding the classification of ergodic line ensembles satisfying a certain Brownian Gibbs property after a parabolic shift.

\subsection{The Airy line ensemble and the Brownian Gibbs property}\label{sec:line ensemble def}
In order to define the Airy line ensemble we first introduce the concept of a \emph{line ensemble} and the Brownian Gibbs line ensembles. We follow the notations of \cite{corwin2014brownian}.

\begin{definition}\label{maindef}
Let $\Sigma$ be a (possibly infinite) interval of $\Z$, and let $\Lambda$ be an interval of $\R$. 
Consider the set $X$ of continuous functions $f:\Sigma\times \Lambda \rightarrow \R$ endowed with the topology of uniform convergence on compact subsets of $\Sigma\times\Lambda$. Let $\mathcal{C}$ denote the sigma-field  generated by Borel sets in $X$.

A  $\Sigma$-\emph{indexed line ensemble} $\mathcal{L}$ is a random variable defined on a probability space $(\Omega,\mathcal{B},\P)$, taking values in $X$ such that $\mathcal{L}$ is a $(\mathcal{B},\mathcal{C})$-measurable function. Intuitively, $\mathcal{L}$ is a collection of random continuous curves (even though we use the word ``line'' we are referring to continuous curves), indexed by $\Sigma$, each of which maps $\Lambda$ into $\R$. We will often slightly abuse notation and write $\mathcal{L}:\Sigma\times \Lambda \rightarrow \R$, even though it is not $\mathcal{L}$ which is such a function, but rather $\mathcal{L}(\omega)$ for each $\omega \in \Omega$. Furthermore, we write $\mathcal{L}_i:=(\mathcal{L}(\omega))(i,\cdot)$ for the line indexed by $i\in\Sigma$.
\end{definition}

We turn now to formulating the Brownian Gibbs property.
\begin{definition}\label{def: maindefBGP}
Let $\{ x_{1}>\dotsm>x_{k}   \}$ and  $\{y_{1}>\dotsm > y_{k}   \} $ be two sets of real numbers. Let $a,b \in \R$ satisfy $a < b$, and  let $f,g:[a,b] \to \R^*$ (where $\R^*=\R\cup\{-\infty,+\infty\}$) be two given continuous functions that satisfy $f(r)>g(r)$ for all $r\in[a,b]$ as well as the boundary conditions $f(a)>x_{1}$, $f(b)>y_{1}$ and $g(a)<x_{k}$, $g(b)<y_{k}$.

The {\it $(f,g)$-avoiding Brownian line ensemble on the interval $[a,b]$ with entrance data  $(x_{1},\ldots,x_{k})$ and exit data $(y_{1},\ldots, y_{k})$} is a line ensemble $\mathcal{Q}$ with $\Sigma=\{1,\ldots, k\}$, $\Lambda=[a,b]$ and with the law of $\mathcal{Q}$ equal to the law of $k$ independent Brownian bridges with diffusion coefficient 1 $\{B_i:[a,b] \to \R\}_{i=1}^{k}$ from $B_i(a) = x_i$ to $B_i(b) = y_i$ conditioned on the event that $f(r) > B_1(r)>B_2(r)>\cdots >B_k(r) > g(r)$ for all $r \in [a,b]$. Note that any such line ensemble $\mathcal{Q}$ is necessarily non-intersecting.

Now fix an interval $\Sigma\subseteq \Z$ and $\Lambda\subseteq \R$ and let $K=\{k_1,k_1+1,\ldots,k_2-1, k_2\} \subset \Sigma$ and $a,b \in \Lambda$, with $a <b$. Set $f=\mathcal{L}_{k_1-1}$ and $g=\mathcal{L}_{k_2+1}$ with the convention that if $k_1-1\notin\Sigma$ then $f\equiv +\infty$ and likewise if $k_2+1\notin \Sigma$ then $g\equiv -\infty$. Write $D_{K,a,b} = K \times (a,b)$ and $D_{K,a,b}^c = (\Sigma \times \Lambda) \setminus  D_{K,a,b}$. A $\Sigma$-indexed line ensemble $\mathcal{L}:\Sigma\times \Lambda \to \R$ is said to have the {\it Brownian Gibbs property} if
\begin{equation*}
\textrm{Law}\Big(\mathcal{L} \big\vert_{D_{K,a,b}} \textrm{conditional on } \mathcal{L} \big\vert_{D_{K,a,b}^c}\Big) = \textrm{Law}(\mathcal{Q}),
\end{equation*}
where $\mathcal{Q}_{i}=\tilde{\mathcal{Q}}_{i-k_1+1}$ and $\tilde{\mathcal{Q}}$ is the $(f,g)$-avoiding Brownian line ensemble on $[a,b]$ with entrance data $\big(\mathcal{L}_{k_1}(s),\ldots,\mathcal{L}_{k_2}(s)\big)$ and exit data  $\big(\mathcal{L}_{k_1}(t),\ldots,\mathcal{L}_{k_2}(t)\big)$. Note that $\tilde{\mathcal{Q}}$ is introduced because, by definition, any such $(f,g)$-avoiding Brownian line ensemble is indexed from $1$ to $k_2-k_1+1$, but we want $\mathcal{Q}$ to be indexed from $k_1$ to $k_2$.
\end{definition}

\smallskip

\begin{definition}\label{def:the Airy line ensemble}
The \emph{Airy line ensemble} is a $ \N\times \R $ indexed line ensemble which we denote by $ \mathcal{A} $. Given $ I\subset \R$, let $ \A(I) = \{\A(i,t)|i\in \N, t\in I  \} $.   The defining property of $ \A  $ is the following: for all $ I=\{t_1\dotsm,t_n  \} $, and $ n\ge 1 $, as a point process on  $ I\times R $, $ \A(I) $ is a determinantal process whose kernel is the \emph{extended Airy$_2 $ kernel} $\KAirytwoext$ such that
\begin{equation}
\KAirytwoext(s,x;t,y)=
\begin{dcases}
\int_0^{\infty} d\lambda\,e^{-\lambda(s-t)} \Ai(x+\lambda)\Ai(y+\lambda) & \textrm{if}\; s \geq t,\\
-\int_{-\infty}^{0} d\lambda\,e^{-\lambda(s-t)} \Ai(x+\lambda)\Ai(y+\lambda) & \textrm{if} \;s<t,
\end{dcases}\label{eq:KAi2}
\end{equation}
where $\Ai(\cdot)$ is the Airy function.
\end{definition}

It is not a priori clear that there exists a continuous line ensemble satisfying the correlation functions in Definition \ref{def:the Airy line ensemble}. \cite{corwin2014brownian}
showed that the Airy line ensemble exists and satisfies the Brownian Gibbs property after a parabolic shift. To be precise, we have
\begin{theorem}[Theorem 3.1 in \cite{corwin2014brownian}]\label{thm:main theorem in ALE paper}
There is a unique $ \N\times \R $ indexed line ensemble $ \mathcal{A} $ satisfying Definition \ref{def:the Airy line ensemble}. Moreover, the line ensemble $ \L $ given by
\begin{equation}\label{eq:parabolic-shift}
\L_i(x) =\frac{1}{\sqrt{2}}(\A_i(x)-x^2 )
\end{equation}
satisfies the  Brownian Gibbs property.
\end{theorem}

\begin{definition}\label{def:Airy$_2$ process}
The top line of the Airy line ensemble is called \emph{Airy$_2$ process}.
\end{definition}

\subsection{Main result and motivation}\label{sec:main result}
%
From the definition of the extended Airy$_2  $ kernel \eqref{eq:KAi2}, the Airy line ensemble and its marginal, the Airy$_2 $ process (i.e. $\mathcal{A}_1$), are  invariant under horizontal shift in the $ x-$coordinate (i.e. are stationary). It is proved  in \cite{prahofer2002scale} that the Airy$_2  $ process is ergodic, and  in fact, satisfies the strong mixing condition. (Technical details regarding ergodicity are provided in Section \ref{sec: Ergodicity}). In this paper, we show that:
\begin{theorem}\label{thm: Main theorem}
The Airy line ensemble is ergodic with respect to horizontal shifts.
\end{theorem}

Theorem \ref{thm: Main theorem} can be seen as a multi-line extension of the ergodicity of the Airy$ _2 $ process. Besides its independent interest, one motivation for this result is our desire to classify stationary ergodic line ensembles which display the Brownian Gibbs property after the parabolic shift given by \eqref{eq:parabolic-shift}.

According to Theorem \ref{thm:main theorem in ALE paper}, $ \mathcal{L}_i(x)=\frac{1}{\sqrt{2}}(\A_i(x)-x^2 ) $ satisfies the Brownian Gibbs property. In \cite{corwin2014brownian}, the authors formulated a conjecture which was originally suggested by Scott Sheffield:
\begin{conjecture}[Conjecture 3.2 of \cite{corwin2014brownian}]\label{conj}
Let $ \Theta=\{\theta_s|s\in \R  \} $ denote the horizontal shift group of  $ \N\times \R $-indexed line ensembles. More precisely, given a line ensemble $ \L $, let
\begin{equation*}
\theta_s\L_i(x)= \L_i(x+s) \;\forall i\in \N,x\in \R.
\end{equation*}
We say that an $\N\times \R$-indexed line ensemble $\mathcal{L}$ is horizontal shift-invariant if $\theta_s\mathcal{L}$ is equal in distribution to $\mathcal{L}$ for each $s \in \R$.
Let $  \mathscr{G}_{\Theta} $ be the set of  Brownian Gibbs line ensemble measure $ \mathcal{L} $'s such that $  2^{1/2} \mathcal{L}_i(x)+x^2  $ is  horizontal shift-invariant. Then as a convex set, the extremal points of $ \mathscr{G}_{\Theta} $ consists of $  \{\mathcal{L}^c| c\in \R \}$ where
\begin{equation}\label{eq:shiftAiry}
\mathcal{L}^c_i(x)=\frac{1}{\sqrt{2}}(\A_i(x)-x^2 ) +c
\end{equation}
  and $\A $ is the Airy line ensemble.
\end{conjecture}

Beyond its intrinsic interest, this conjecture is worth investigating in light of its possible use
as an invariance principle for deriving convergence of systems to the Airy line ensemble. As such,
the characterization could serve as a route to universality results. For example, Section 2.3.3 in \cite{corwin2013kpz} suggests a possible route to prove that the KPZ line ensemble converges to the Airy line ensemble minus a parabola (and hence the narrow-wedge initial data KPZ equation converges to the Airy$ _2 $ process minus a parabola) based on the above conjecture.

\bigskip
In Section \ref{sec:Extremal} we will prove  that Theorem \ref{thm: Main theorem} implies:
\begin{theorem}\label{thm:extremal}
$  \mathcal{L}^c  $ defined in Conjecture \ref{conj} are extremal Brownian Gibbs line ensembles.
\end{theorem}
Theorem \ref{thm:extremal} reduces Conjecture \ref{conj} to the following:
\begin{conjecture}\label{conj:reduced one}
Given $ c\in \R $, there is a unique Brownian Gibbs line ensemble $ \mathcal{L} \in \mathscr{G}_{\Theta} $ such that
\[
 \E\big[ \mathcal{L}_1(x)+ 2^{-1/2} x^2 \big]=c
\]for all $ x\in \R $.
\end{conjecture}

When Conjecture \ref{conj} was formulated in \cite{corwin2014brownian}, it was not shown that $ \mathcal{L}^c $ defined in \eqref{eq:shiftAiry} are extremal.  Therefore even if the uniqueness in Conjecture \ref{conj:reduced one} was established, it would not necessarily follow that extremal points in $ \mathscr{G}_{\Theta} $ were related to the Airy line ensemble. Theorem \ref{thm:extremal} rules out the possibility that the Airy line ensemble is a nontrivial convex combination of extremal points in $\mathscr{G}_{\Theta}$, thus reducing Conjecture \ref{conj} to Conjecture \ref{conj:reduced one}.

\subsection{Ergodicity of the Airy line ensemble}\label{sec: Ergodicity}
We recall some basic facts in ergodic theory in the context of $ \N\times \R $-indexed line ensemble.  As a matter of convention, all  line ensembles are assumed to be indexed by $ \N\times \R  $ unless otherwise noted.\\

Recall $ \Theta=\{\theta_s|s\in \R  \} $ is the horizontal shift group of  $ \N\times \R $-indexed line ensembles.

\begin{definition}\label{def:ergodicity}
Suppose $ \L $ is a horizontal shift-invariant line ensemble on the probability space $ (\Omega, \mathcal{C} , \P ) $. We say  $ A\in \mathcal{C}$ is shift-invariant if $ \theta_sA:=\{\theta_s\omega\,| \omega \in A  \}  =A $ for all $ s \in\R $. Then $ \L $ is ergodic if for all shift-invariant $ A $, $ \P[A] =0 $ or 1.
\end{definition}
It is well known that ergodicity follows from the  \emph{strong mixing condition}.
\begin{definition}\label{def:mixing condition}
Suppose $ \L $ is a horizontal shift-invariant line ensemble on the probability space $ (\Omega, \mathcal{C} , \P ) $.  $ \L  $ is said to satisfy the strong mixing condition if for all $ A ,B\in \mathcal{C} $,
\[
\lim\limits_{T\to \infty }\P[ \theta_TA,B  ]=\P[A]\P[B].
\]
\end{definition}
\begin{proposition}\label{thm:mixing implies ergodic}
If a horizontal shift-invariant line ensemble $ \L $ satisfies the strong mixing condition, it is ergodic.
\end{proposition}
\begin{proof}
Suppose $ A $ is shift-invariant. Then $ \P[ \theta_TA,A  ]=\P[A] $. Let $ T $ tend to infinity, by the strong mixing condition, $ \P[A]=\P[A]^2 $, which means $ \P[A] =0$ or 1.
\end{proof}

To prove Theorem \ref{thm: Main theorem}, we actually prove a stronger result:
\begin{proposition}\label{thm:strong mixing condition}
The Airy line ensemble satisfies the strong mixing condition given in Definition \ref{def:mixing condition}.
\end{proposition}

Now we consider the Airy line ensemble $ \mathcal{A} $. Fix $ m\in \N $ and  $ t_1<t_2<\cdots< t_m  $. The Airy line ensemble restricted to $  \N\times \{t_1, t_2,\cdots, t_m \}$ is a point process on $ \{t_1, t_2,\cdots, t_m  \}\times \R $. For $ 1\le i\le m $ and $ k_i\in \N $,  suppose  $\mathcal{I} = \{I^j_i|1\le j \le k_i \} $ is a collection of intervals on $ \N\times \R $ satisfying
\begin{align}\label{eq:lower bound for interval}
&\{I^j_i \}_{1\le j \le k_i} \textrm{ are  disjoint intervals on }  \{t_i \}\times \R \textrm{ for all } i, \textrm{ and} \nonumber\\
&\inf\Big\{x \big| \exists 1\le i\le m, 1\le j \le k_i  \textrm{ such that } (t_i,x) \in I_i^j \Big \}>-\infty.
\end{align}

Let $ N^{I_i^j} $ denote the number of particles in  $ I^j_i $ and let
\[
\mathcal{G}(t_1,t_2\dotsm,t_m,\mathcal{I})=\Big\{ A\in \mathcal{C}| A\in  \sigma\Big( \{N^{I_i^j}|  1\le i\le m,1\le j\le k_i \}\Big) \Big \}.
\]

We prove Proposition \ref{thm:strong mixing condition} by showing that
\begin{lemma}\label{lemma: mixing for particles number in interval}
Fix $ m\in \N $, $ t_1<t_2<\cdots< t_m  $ and $ \mathcal{I}$,  for $ A,B \in \mathcal{G}(t_1,t_2,\dotsm,t_m,\mathcal{I} ) $, 
\begin{equation}\label{eq:mixing for intervals}
\lim\limits_{T\to \infty }\P[  \theta_T (A),B  ]=\P[A]\P[B].
\end{equation}
\end{lemma}
\begin{proof}[Proof of Proposition \ref{thm:strong mixing condition} based on  Lemma \ref{lemma: mixing for particles number in interval} ]
Recall a result in measure theory (see, for example,\cite{halmos1950measure}),
\begin{theorem}
Let $ (X,\mathcal{C},\P)  $ be a probability space. Let $ \mathcal{B}\subset \mathcal{C} $ an algebra generating $ \mathcal{C} $.
Then for all $ A \in \mathcal{C} $and $ \eps >0 $, we can find $ A'\in \mathcal{B} $ such that $ \mu(A\Delta A')\le\eps $. Here $A\Delta A'=A\backslash A'+A'\backslash A$ is the symmetric difference of $ A $ and $ A' $.
\end{theorem}

Let $\mathcal{F}$ be the $\sigma$-algebra generated by the Airy line ensemble and
$ \mathcal{F}'$ be the union of all $\mathcal{G}(t_1,t_2\dotsm,t_m,\mathcal{I})$
where $(t_1,t_2\dotsm,t_m)$ varies over all finite collections of real numbers and $\mathcal I$ varies over all finite collections of intervals satisfying \eqref{eq:lower bound for interval}.
Since $$\mathcal{G}(t_1,t_2\dotsm,t_m,\mathcal{I})\cup\mathcal{G}(t'_1,t'_2\dotsm,t'_n,\mathcal{I'})\subset \mathcal{G}(t_1,t_2\dotsm,t_m,t'_1,t'_2\dotsm,t'_n, \mathcal{I}\cup\mathcal I')$$ for all choices of $(t_1,t_2\dotsm,t_m,\mathcal{I})$ and $(t'_1,t'_2\dotsm,t'_n,\mathcal{I'})$, $\mathcal F'$ is an algebra. Furthermore, if we restrict the Airy line ensemble to $\N\times\Q$, then $\forall i\in \N, r\in\Q$, $\mathcal A_i(r)$ is measurable with respect to the $\sigma$-algebra generated by $\mathcal F'$. By the continuity of the Airy line ensemble, $\mathcal F$ equals to the $\sigma$-algebra generated by $\mathcal F'$.

Therefore for $ A,B\in \mathcal{F}$ and $\eps>0$, there exist $A',B'\in \mathcal F'$ such that $ \P[A\Delta A']<\eps $, $  \P[B\Delta B']<\eps $.
By Lemma \ref{lemma: mixing for particles number in interval},
\begin{equation*}
\lim\limits_{T\to \infty }\P[  \theta_T (A'),B'  ]=\P[A']\P[B'].
\end{equation*}
Since $ \eps $ can be arbitrarily small,
\begin{equation*}
\lim\limits_{T\to \infty }\P[  \theta_T (A),B   ]=\P[A ]\P[B ].
\end{equation*}
\end{proof}
\subsection{Strategy and outline}
Now we sketch our strategy for proving our main results, Theorem \ref{thm: Main theorem} and \ref{thm:extremal}.
From the argument above, Lemma \ref{lemma: mixing for particles number in interval} implies Proposition \ref{thm:strong mixing condition}, which further implies Theorem \ref{thm: Main theorem}. Thus Theorem \ref{thm: Main theorem} boils down to proving Lemma \ref{lemma: mixing for particles number in interval}, which is the content of Section \ref{sec: Mixing}.

Since $ N^{I_i^j} $ are discrete random variables, their joint distribution is governed by the moment generating function.
By a standard fact in determinantal point process (Lemma \ref{lem:generating function}), the generating function can be expressed as a Fredholm determinant of $ \KAirytwoext $ on $L^2 (\{  t_1,t_2,\dotsm,t_m \}\times \R ) $.

Since there is a time shift $ T $ in Lemma \ref{lemma: mixing for particles number in interval}, we need to consider $ 2m $ time moments. The trace class operator involved in the left hand side of \eqref{eq:mixing for intervals}  can be regarded as a $ 2m $ by $ 2m $ operator valued matrix. Based on the known estimates of Airy$_2$ kernel, one can show that as a 2 by 2 block matrix of block size $ m $, the ``off diagonal'' terms  will vanish as $ T\to\infty $. This shows the factorization of the generating function (Lemma \ref{lem: splitting}). Then using complex analysis of several variables one can extract Lemma \ref{lemma: mixing for particles number in interval} from Lemma \ref{lem: splitting}.

To prove Theorem \ref{thm:extremal}, we follow the standard method in \cite{georgii2011gibbs}. As shown in Chapter 14 of \cite{georgii2011gibbs}, given a translation group and a Gibbs specification,  for the set of translation invariant Gibbs measures defined on $ \mathbb{S}=\Z^d $, extremal points coincide with ergodic ones with respect to the translation group. In our case, the underlying space is $ \mathbb{S}= \R\times \N  $. Thanks to the fact that the Airy line ensemble is a collection of \emph{continuous} curves \cite{corwin2014brownian}, the same argument in \cite{georgii2011gibbs} works for the Airy line ensemble.  The details are provided in Section \ref{sec:Extremal}.

\subsection{Acknowledgments}
We thank Alan Hammond, Scott Sheffield and Herbert Spohn for discussions. IC was partially supported by the NSF through DMS-1208998, by Microsoft Research through the Schramm Memorial Fellowship, and by the Clay Mathematics Institute through a Clay Research Fellowship.


\section{Strong mixing for $ N^{I^j_i} $}\label{sec: Mixing}
\subsection{Express the generating function by Fredholm determinant}\label{sec:Fredholm Determinant}
We first recall a useful formula from determinantal point process. It appears in \cite{borodin2009determinantal} as formula (2.4).

\begin{lemma}\label{lem:generating function}
Suppose $ X $ is a determinantal point process on a locally compact space of $ \mathcal{X} $ with kernel $ K $. Let $ \phi  $ be a function on $ \mathcal{X} $ such that the kernel $ (1-\phi(x)) K(x,y) $ defines a trace class operator $ (1-\phi)K $  in $ L^2(\mathcal{X}) $. Then
\begin{equation}
\E\Big[\prod\limits_{x_i\in X} \phi(x_i)\Big] =\det\big(  \1  -(1-\phi)K \big).
\end{equation}
\end{lemma}
Fix  $ z_1,z_2,\cdots, z_m \in \C$ such that $ |z_i|\le 1 $ for all $ i $. Let  $I_1,I_2, \cdots, I_m $ be a family of pairwise disjoint subsets of $ \mathcal{X} $ and $ Q $ is a multiplication operator defined by
\begin{equation}
Qf(x)=\sum\limits_{i=1}^{m}(1-z_i )\1_{x\in I_i} f(x).
\end{equation}
Denote $ N_{I_i} (1\le i \le m)  $  to be the number of particles in $ I_i $ of  the random configuration $ X $. Specifying
\[\phi(x)=\sum\limits_{i=1}^{m} z_i\1_{x\in I_i}  + \1_{x\in \cap_{i=1}^m I_i^c}   \]
 in Lemma \ref{lem:generating function} leads to
\begin{equation}\label{eq:generating function}
\E\Big[\prod\limits_{i=1}^m z_i^{N_{I_i}}  \Big] = \det(I-QK)_{L^2(\mathcal{X})}.
\end{equation}

\subsection{Proof of Lemma \ref{lemma: mixing for particles number in interval} }\label{sec:Spitting}
Let the  \emph{Airy Hamiltonian}  be defined
as \[H=-\Delta +x.\] $H$ has the shifted Airy functions $\Ai_\lambda(x)=\Ai(x-\lambda)$
as its generalized eigenfunctions: $H\!\Ai_\lambda(x) = \lambda\!\Ai_\lambda(x)$. Define
the \textit{Airy$_2$ kernel} $\KAirytwo$ as the projection of $H$ onto its negative generalized
eigenspace:
\[\KAirytwo(x,y)=\int_0^{\infty} d\lambda \Ai(x+\lambda)\Ai(y+\lambda).\]

Consider $ t_1<t_2<\cdots< t_n  $, $ t_i\in \R $. 
For $ 1\le i\le n $,  suppose  $\{I^j_i \}_{1\le j \le k_i}   $ are intervals  on $ \{t_i \}\times \R $ satisfying the condition in \eqref{eq:lower bound for interval} and
\begin{equation}\label{eq:lower bound}
M_0=-\inf\{x\,| (t_i,x) \in I_i^j \textrm{ for some } i,1\le j\le k_i  \}.
\end{equation}
Let $ N^{I_i^j} $ be the number of particles in the interval $ I^j_i $.

From Lemma \ref{lem:generating function} and \eqref{eq:generating function}, for $ \{ z_i^j  | 1\le i \le n,1 \le j\le k_i, |z_i^j |\le 1   \} $
\begin{equation}\label{eq:number of particle}
\E\Big[\prod\limits_{i=1}^n\prod\limits_{j=1}^{ k_i} (z_i^j)^{N^{I^j_i}}\Big] = \det ( I-Q\KAirytwoext)_{L^2(  \{t_1, t_2,\cdots, t_n  \} \times \R )},
\end{equation}
where $ Q $ is a multiplication operator defined as
\begin{equation}\label{eq: def of Q}
Qf(t,x)=\sum_{\substack{
1\le i \le n\\
1\le j\le k_i
}}
(1-z^j_i )\1_{(t,x)\in I^j_i} f(t,x)
\end{equation}
for all $ f\in  L^2(  \{t_1, t_2,\cdots, t_n  \} \times \R )$.\\

For the rest of the section we regard $ Q\KAirytwoext $ as an $ n\times n $ operator valued matrix where
\begin{equation}\label{eq: matrix entry}
 [Q\KAirytwoext]_{ij}= Q_{t_i} e^{(t_i-t_j)H} \KAirytwo \1_{i\ge j}+Q_{t_i} e^{(t_i-t_j)H} (\KAirytwo-I)  \1_{i < j}.
\end{equation}
We write $\det\big(I-[Q\KAirytwoext]_{1\le i,j\le n}\big)$ to be the Fredholm determinant of $Q\KAirytwoext$. For $1\le m\le n$, we write $[Q\KAirytwoext]_{1\le i,j\le m}$ as the operator corresponding to the submatrix of $Q\KAirytwoext$ which consists of  entries indexed by $\{1,\dotsm,m\}\times \{1,\dotsm,m\}$ and $\det(I - [Q\KAirytwoext]_{1\le i,j\le m})$ as its Fredholm determinant. $\det\big(I-[Q\KAirytwoext]_{m+1\le i,j\le n}\big)$ is defined in the same manner.

According to the convention above,
\begin{equation}\label{eq: number of particle}
\E\Bigg[\prod\limits_{i=1}^n\prod\limits_{j=1}^{ k_i} (z_i^j)^{N^{I^j_i}}  \Bigg]=\det\big(I-[Q\KAirytwoext]_{1\le i,j\le n}\big).
\end{equation}

Now we assume $ n=2m $ is an even number and $ \{t_1< t_2<\cdots< t_m \} $ is a fixed set of real numbers.  $ \{t_{m+1},t_{m+2}\dotsm, t_{n}  \} $ is a shift of $ \{t_1, t_2,\cdots, t_m \} $ by $ T $. To be precise,  $ t_{m+i}=t_i +T$ for all $ 1\le i\le m $.

Similar to \eqref{eq: number of particle}, we have
\begin{align}\label{eq: number of particle 2}
\E\Bigg[\prod\limits_{i=1}^m\prod\limits_{j=1}^{k_i} (z_i^j)^{N^{I^j_i}} \Bigg ]&=\det\big(I-[Q\KAirytwoext]_{1\le i,j\le m}\big),\nonumber\\
\E\Bigg[\prod\limits_{i=m+1}^n\prod\limits_{j=1}^{k_i} (z_i^j)^{N^{I^j_i}} \Bigg]&=\det\big(I-[Q\KAirytwoext]_{m+1\le i,j\le n}\big).
\end{align}

\begin{lemma}\label{lem: splitting}
Let
\begin{align}\label{eq: splitting}
R(z,T)&=\E\Bigg[\prod\limits_{i=1}^n\prod\limits_{j=1}^{k_i} (z_i^j)^{N^{I^j_i}} \Bigg]-\E\Bigg[\prod\limits_{i=1}^m\prod\limits_{j=1}^{k_i} (z_i^j)^{N^{I^j_i}} \Bigg ]\E\Bigg[\prod\limits_{i=m+1}^n\prod\limits_{j=1}^{k_i} (z_i^j)^{N^{I^j_i}} \Bigg]\nonumber\\
&=\det\big(I-[Q\KAirytwoext]_{1\le i,j\le n}\big)-\det\big(I-[Q\KAirytwoext]_{1\le i,j\le m}\big)\det\big(I-[Q\KAirytwoext]_{m+1\le i,j\le n}\big).
\end{align}
then $\lim\limits_{T\to\infty} R(z,T)=0 $ for all $z$.
\end{lemma}

We postpone the proof of Lemma \ref{lem: splitting} to Section  \ref{sec:proof of splitting}. Now we prove Lemma \ref{lemma: mixing for particles number in interval} based on it.
Let $ N_i^j\in \N $ for all $ 1\le i\le n, 1\le j\le k_i  $.  Define the partial differential operators
\begin{equation}
\partial_{1,n}=\frac{\partial^{ \sum\limits_{i=1}^n\sum\limits_{j=1}^{k_i} N_i^j } }{ \prod\limits_{i=1}^n\prod\limits_{j=1}^{k_i} \partial (z_i^j)^{N_i^j} },\quad
\partial_{1,m}=\frac{\partial^{ \sum\limits_{i=1}^m\sum\limits_{j=1}^{k_i} N_i^j } }{ \prod\limits_{i=1}^m\prod\limits_{j=1}^{k_i} \partial (z_i^j)^{N_i^j} },\quad
\partial_{m+1,n}=\frac{\partial^{ \sum\limits_{i=m+1}^n\sum\limits_{j=1}^{k_i} N_i^j } }{ \prod\limits_{i=m+1}^n\prod\limits_{j=1}^{k_i} \partial (z_i^j)^{N_i^j} }.
\end{equation}
Then  by Lemma \ref{lem: splitting}
\begin{align}\label{eq:partial derivative}
&\partial_{1,n}\E\Bigg[\prod\limits_{i=1}^n\prod\limits_{j=1}^{k_i} (z_i^j)^{N^{I^j_i}}  \Bigg] \nonumber\\
=&\partial_{1,m}\E\Bigg[\prod\limits_{i=1}^{m}\prod\limits_{j=1}^{k_i} (z_i^j)^{N^{I^j_i}} \Bigg ]\partial_{m+1,n}\E\Bigg[\prod\limits_{i=m+1}^n\prod\limits_{j=1}^{k_i} (z_i^j)^{N^{I^j_i}} \Bigg ]
+\partial_{1,n} R(z,T).
\end{align}

 $ R(z,T) $ is analytic function bounded by 2 when $ |z|< 1 $.
 By Montel's theorem \cite[Theorem 1.4.31]{severalcomplex}, every subsequence of $ R(z,T) $ has a further  subsequence which converges locally uniformly. This implies that $\lim\limits_{T\to \infty} R(z,T) =0$ locally uniformly in $ |z|< 1 $.
By the Cauchy inequalities \cite[Theorem 1.3.3]{severalcomplex}, the magnitude of partial derivatives of an analytic function at a certain point is controlled by the magnitude of the  function around that point. Therefore $ \partial_{1,n} R(z,T) $ converges uniformly for $ |z|<1 $. In particular, $ \lim\limits_{T\to\infty}\partial_{1,n} R(0,T) =0$.

Let
\begin{align*}
A'&=\{ N^{I_i^j} =N_i^j \textrm{ for all } 1\le i\le m, 1\le j\le k_i \},\\
B'&=\{ N^{I_i^j} =N_{i+m}^j \textrm{ for all } 1\le i\le m, 1\le j\le k_i \}.
\end{align*}
Taking $ z=0 $ in  \eqref{eq:partial derivative} we have
\begin{align}\label{eq: split 1}
\P[ A', \theta_TB']=\P[A']\P[\theta_TB' ]+o(1)=\P[A']\P[B' ]+o(1)  \textrm{ as }T\to\infty,
\end{align}where the second equality is due to the horizontal shift invariance of the Airy line ensemble.

Therefore for all $ A, B \in \mathcal{G}(t_1,\dotsm,t_m,\mathcal{I}) $,
\begin{align}\label{eq: split 2}
\lim\limits_{T\to\infty }\P[ A,\theta_TB]= \P[A]\P[B ].
\end{align}
This finishes the proof of Lemma \ref{lemma: mixing for particles number in interval}.

\subsection{Proof of Lemma  \ref{lem: splitting}}\label{sec:proof of splitting}

The proof of Lemma \ref{lem: splitting} uses the following fact.
\begin{lemma}\label{lem:tracedecay}
Let $ P_a $ be the multiplication operator of $ \1_{x>a} $ and $ y>0 $. Then $ P_a\KAirytwo $, $ P_ae^{-yH}(I-\KAirytwo) $ and $ P_ae^{yH} \KAirytwo $ are trace class operators. Moreover,
\begin{equation}\label{eq: tracedecay}
\lim\limits_{y\to\infty} \| P_a e^{-yH}(I-\KAirytwo)\|_1 =\lim\limits_{y\to\infty} \|P_a e^{yH} \KAirytwo\|_1  =0.
\end{equation}
\end{lemma}
\begin{proof}
Let $ \varphi(x) = (1+x^2) ^{1/2} $ and define the multiplication operator $ Mf(x)= \varphi (x) f(x) $.

In the proof of Proposition 3.2 in \cite{corwin2013continuum}, it was shown that
 \begin{equation}
  \Big \|P_ae^{-yH}M \Big \|_2 \le C.
\end{equation}where $  \|\cdot \|_2$ is the Hilbert-Schmidt norm.
Actually by taking $ n=2 $ in the proof of Proposition 3.2 in \cite{corwin2013continuum}, we have
 \[
\Big \|P_ae^{-yH}M \Big \|_2=   \Big \| \big[e^{-yH}-(I-P_a)e^{-yH}  \big]M \Big \|_2 \le C.
  \]
  On the other hand
 \begin{equation}
    \begin{split}
      \|M^{-1}e^{y H}\KAirytwo\big\|^2_2&=\int_{\R^2}dx\,dz\int_{(-\infty,0]^2}d\lambda\,d\tilde\lambda\,
      \varphi(x)^{-2}e^{(\lambda+\tilde\lambda)y}\Ai(x-\lambda)\Ai(z-\lambda)\\
      &\hspace{3in}\cdot\Ai(x-\tilde\lambda)\Ai(z-\tilde\lambda)\\
      &=\int_{-\infty}^\infty dx\int_{-\infty}^0d\lambda\,\varphi(x)^{-2}e^{2\lambda
        y}\Ai(x-\lambda)^2\\&\le  \frac{c}{2y}\|\varphi^{-1}\|_2^2,
    \end{split}
  \end{equation}where $c=\max_{x\in\R}\Ai(x)^2<\infty$.

Similarly,
\begin{equation}
\|M^{-1}e^{-y H}(I-\KAirytwo)  \big\|_2^2 \le  \frac{c}{2y}\|\varphi^{-1}\|_2^2.
\end{equation}

Since $ \|AB\|_1\le   \|A\|_2 \|B\|_2 $,
\[\| P_a K \|_1\le \| P_ae^{yH} M \|_2  \| M^{-1}e^{-yH} K \|_2 <\infty. \]
Fix $ b>0 $, there exists a constant $ C_b $ such that
\begin{align*}
&\| P_a e^{yH}  K \|_1\le \| P_ae^{-bH} M \|_2  \| M^{-1}e^{-(y-b)H} K \|_2 \le \frac{C_b}{y},\\
&\| P_a e^{-yH}(I-K)   \|_1\le \| P_ae^{-bH} M \|_2  \| M^{-1}e^{-(y-b)H}(I-K) \|_2 \le \frac{C_b}{y}.
\end{align*}
This finishes the proof.
\end{proof}

To prove Lemma \ref{lem: splitting}, we come back to \eqref{eq: matrix entry}.

\begin{equation}
R(z,T)=\det\big(I-[Q\KAirytwoext]_{1\le i,j\le n}\big)-\det\big(I-[Q\KAirytwoext]_{1\le i,j\le m}\big)\det\big(I-[Q\KAirytwoext]_{m+1\le i,j\le n}\big).
\end{equation}

By condition \eqref{eq:lower bound for interval}, we can always replace $ Q_{t_i}$ by $Q_{t_i}P_a $ for $ a<-M_0 $. Therefore from Lemma \ref{lem:tracedecay},
\begin{align*}
&\lim\limits_{T\to \infty}\| Q_{t_i} e^{(t_i-t_j)H} \KAirytwo\|_1 =0,\;& 1\le j\le m<i\le n,\\
&\lim\limits_{T\to \infty}\| Q_{t_i} e^{(t_i-t_j)H} (\KAirytwo-I)\|_1 =0,&1\le i\le m<j\le n.\\
\end{align*}

Since Fredholm determinant is continuous respect to trace norm, $ \lim\limits_{T\to\infty} R(z,T) =0 $.



\section{Extremal Gibbs measure and ergodicity}\label{sec:Extremal}
\subsection{The Gibbs property of the Airy line ensemble}\label{sec: Gibbs property of ALE}
We first adjust the notations to make it consistent with those in \cite{georgii2011gibbs}.

Let $ S=\N\times\R $ denote the parameter set and $ (\Omega,\mathscr{F},\P  ) $ denote the probability space of the Airy line ensemble. Canonically, we choose $ \Omega $ to be the set of continuous functions from $ S$ to $ \R $ and $ \mathscr{F} $ be the Borel $ \sigma-$algebra with respect to the locally uniformly convergence topology. Let $\mathscr{L} $ be the set of all allowable finite non-empty subsets of $ S $. Here allowable sets are those of the form $\{ k_1,\dotsm,k_2\} \times (a,b) $. Let $ \mathscr{T}_{\Lambda}=\mathscr{F}_{S\backslash\Lambda} $ where  $ \Lambda $ runs through $ \mathscr{L} $ and $\mathscr{F}_{S\backslash\Lambda}  $ be the $ \sigma$-algebra of the Airy line ensemble restricted to $ S\backslash\Lambda$.

From Theorem \ref{thm:main theorem in ALE paper}, the Airy line ensemble satisfies certain Gibbs property, which we formulate now.

Let $\{ x_{1}>\dotsm>x_{k}   \}$ and  $\{y_{1}>\dotsm > y_{k}   \} $ be two sets of real numbers. Let $a,b \in \R$ satisfy $a < b$, and  let $f,g:[a,b] \to \R^*$ (where $\R^*=\R\cup\{-\infty,+\infty\}$) be two given continuous functions that satisfy $f(r)>g(r)$ for all $r\in[a,b]$ as well as the boundary conditions $f(a)>x_{1}$, $f(b)>y_{1}$ and $g(a)<x_{k}$, $g(b)<y_{k}$. The \textit {shifted $(f,g)$-avoiding Brownian line ensemble on the interval $[a,b]$ with entrance data  $(x_{1},\ldots,x_{k})$ and exit data $(y_{1},\ldots, y_{k})$} is a line ensemble $ \mathcal{L} $ such that $ 2^{-1/2}(\L-x^2)  $ is a $\big(2^{-1/2}(f-x^2),2^{-1/2}(g-x^2)\big)$-avoiding Brownian line ensemble on the interval $[a,b]$ with entrance data  $\big(2^{-1/2}(x_1-a^2),\ldots,2^{-1/2}(x_k-a^2)\big)$ and exit data $\big(2^{-1/2}(y_1-b^2),\ldots,2^{-1/2}(y_k-b^2)\big)$.

\begin{definition}\label{def: Gibbs specilization}
Suppose $ \Lambda=\{ k_1,\dotsm,\kappa_2\} \times (a,b) $ and $ \gamma_\Lambda $ is a probability kernel from $ (\Omega,\mathscr{T}_{\Lambda}) $ to $ (\Omega,\mathscr{F} ) $ defined as follows:

For $ \omega\in \Omega $, $ \gamma_\Lambda(\cdot| \omega ) $ coincides with $ \omega $ outside $ \Lambda $; in $ \Lambda $, $ \gamma_\Lambda(\cdot| \omega ) $ is the law of the \emph{shifted} $ \big(\omega(k_1-1, \cdot)\big|_{[a,b]}, \omega(k_2+1, \cdot)\big|_{[a,b]}  \big)   $-avoiding Brownian line ensemble on $[a,b]$ with entrance data $\big(\omega(k_1,a)\ldots,\omega(k_2,a)\big)$ and exit data  $\big(\omega(k_1,b)\ldots,\omega(k_2,b)\big)$. Here we make the convention that $ \omega(0,x)=\infty $.
\end{definition}
$ \gamma=\big( \gamma_\Lambda \big)_{\Lambda\in \mathscr{L}}$ is the family of Gibbs specifications that describes the Gibbs property of the Airy line ensemble. Since the Airy line ensemble is stationary, $ \gamma $ is also horizontal shift-invariant, which means
\begin{equation}\label{eq: shift-invariant of gamma}
\gamma_{\Lambda+T} (\theta_TA |\theta_T\omega  )= \gamma_{\Lambda} (A |\omega  ), \;\;\;( \Lambda\in \mathscr{L},T\in \R, \omega\in \Omega ).
\end{equation}

With these notations, Theorem \ref{thm:extremal} can be formulated in this way.
\begin{theorem}\label{thm: extremal'}
Let $ \mathscr{G}_{\Theta}(\gamma)  $ be the simplex  of all probability measures $ \mu $ on $ (\Omega,\mathscr{F}) $ such that
\begin{equation}
\mu(\theta_TA)=\mu(A) \textrm{ and } \mu(A|\mathscr{T}_\Lambda ) =\gamma(A|\cdot)  \quad \mu \textrm{ a.s. for all } A\in \mathscr{F}, \Lambda\in \mathscr{L}\textrm{ and } T \in \R .
\end{equation}
Then the Airy line ensemble is an extreme  point of $  \mathscr{G}_{\Theta}(\gamma)  $.
\end{theorem}

\subsection{Proof of Theorem \ref{thm: extremal'}}\label{sec: proof of extremal}

Theorem \ref{thm: extremal'} follows from a standard result in ergodic theory which is proved as Corollary 7.4 in \cite{georgii2011gibbs}.

\begin{lemma}\label{lem: tiviality}
Let $ (\Omega,\mathscr{F}) $ be a measurable space, $ \Pi $ a non-empty set of probability kernels from $ \mathscr{F} $ to $ \mathscr{F} $, and
  \begin{equation}
\mathscr{P}_{\Pi} =\big\{ \mu\in \mathscr{P}(\Omega,\mathscr{F} ): \mu\pi=\mu \textrm{ for all } \pi\in \Pi \big \},
\end{equation}
be the convex set of all $ \Pi-$invariant probability measures on $ (\Omega,\mathscr{F} ) $. Let $ \mu\in \mathscr{P}_\Pi $ be given and $ \mathscr{I}_\Pi(\mu)=\displaystyle{\cap_{\pi\in\Pi}} \mathscr{I}_\pi(\mu) $ where $ \mathscr{I}_\pi(\mu)=\big\{  A\in \mathscr{F}: \pi(A|\cdot ) = \1_{A} \; \mu\textrm{-a.s.}\big\} $.
Then $ \mu $ is extreme if and only if $ \mu $ is trivial on $ \mathscr{I}_\Pi(\mu) $.
\end{lemma}

Define a family  $ \hat{\Theta} =  \{ \hat{\theta}_t: t\in \R  \} $ of  probability kernels $ \hat{\theta}_T $ from $ \mathscr{F} $ to $ \mathscr{F} $ by
\[  \hat{\theta}_t(A|\omega) = \1_A(\theta_t\omega)\quad \quad (t\in \R, A\in \mathscr{F}, \omega\in \Omega ).  \] Using the notation of Lemma  \ref{lem: tiviality},  $ \mathscr{P}_{\hat{\Theta}}  $  is the set  of all horizontal shift-invariant measures.

By definition, a probability measure $ \mu$ belongs to $ \mathscr{G}_\Theta(\gamma) $ if and only if $ \mu $ is preserved by all probability kernels in
\[ \Pi= \{\gamma_\Lambda:\gamma\in \mathscr{L}\}\cup \mathscr{P}_{\hat{\Theta}} . \] Therefore $ \mu $ is an extreme in $ \mathscr{G}_\Theta(\gamma) $ if and only if $ \mu $ is trivial on $ \mathscr{I}_\Pi(\mu) $.

We claim that  $ \P  $ (the probability measure for the Airy line ensemble) is trivial on $ \cap_{t\in \Q } \mathscr{I}_{\hat{\theta}_t }(\mu) $, which is sufficient to prove Theorem \ref{thm: extremal'}.

To show the triviality, for given $ A\in \cap_{t\in \Q } \mathscr{I}_{\hat{\theta}_t }(\mu) $, let $B= \cup_{t\in \Q}\theta_t A $. Then $ \theta_tB=B $ for all $ t\in \Q $ and
$ \P[A\Delta B]\le \sum\limits_{t\in\Q} \P[A\Delta\theta_tA ] =0$. Therefore $ \P[A]=\P[B] $.

As in the proof of Proposition \ref{thm:mixing implies ergodic}, $ \P[B]=\lim\limits_{n\to\infty}\P[\theta_nB,B ]=\P^2[B]  $, which means $ \P[A]=\P[B]=0 $ or 1.


\bibliographystyle{plain}
\bibliography{ALE}

\filbreak
\begingroup
\small
\parindent=0pt

\bigskip
\vtop{
\hsize=5.3in
Ivan Corwin\\
Department of Mathematics, Columbia University\\
2990 Broadway, New York, NY 10027, USA\\
Clay Mathematics Institute\\
10 Memorial Blvd. Suite 902, Providence, RI 02903, USA\\
Department of Mathematics, Massachusetts Institute of Technology\\
77 Massachusetts Avenue, Cambridge, MA 02139-4307, USA\\
ivan.corwin@gmail.com
}

\bigskip

\vtop{
\hsize=5.3in
Xin Sun\\
Department of Mathematics, Massachusetts Institute of Technology\\
77 Massachusetts Avenue, Cambridge, MA 02139-4307, USA\\
xinsun89@math.mit.edu
}
\endgroup \filbreak
\end{document}